\documentclass{amsart}

\usepackage{amssymb}
\usepackage{color}
\usepackage{amsxtra} 
\usepackage{mathtools}

\newtheorem{theorem}{Theorem}[section]
\newtheorem{lem}[theorem]{Lemma}

\newtheorem*{theorem*}{Theorem}
\newcommand{\RightEqNo}{\let\veqno\@@reqno} 

\usepackage{environ}
\makeatletter
\NewEnviron{Ralign}{\tagsleft@false\begin{align}\BODY\end{align}}
\makeatother

\theoremstyle{definition}

\theoremstyle{remark}

\numberwithin{equation}{section}

\begin{document}

\title[Nonuniqueness of solutions of the Navier-Stokes equations]{Nonuniqueness of solutions of the Navier-Stokes equations on negatively curved Riemannian manifolds}

\author{Leandro A. Lichtenfelz}
\address{Department of Mathematics, University of Notre Dame, South Bend, IN 46556, USA} 
\email{llichte1@nd.edu} 


\date{\today} 


\begin{abstract}
In a well-known work, M. Anderson constructed a Hadamard manifold $(M^n, g)$ which carries non-zero $L^2$ harmonic $p$-forms when $p \neq n/2$, thus disproving the Dodziuk-Singer conjecture. In this paper, we use the manifold $(M^3, g)$ in order to solve another problem in geometric analysis, namely the nonuniqueness of solutions of Leray-Hopf type of the Navier-Stokes equations.
\end{abstract}

\maketitle

\section{Introduction} 
\label{sec:Intro}

The aim of this paper is to prove that the Cauchy problem for the Navier-Stokes equations is ill-posed on the Riemannian manifold $(M^3, g)$ constructed by Anderson \cite{And} in his counterexample to the Dodziuk-Singer conjecture (see also \cite{DonFred}, \cite{Luck}, \cite{Dod}, \cite{Gromov}). This answers a question raised by Chan and Czubak in \cite{CC}, pursued also by Khesin and Misiolek \cite{KM}.

In order to provide the appropriate context, we begin with a brief review of some pertinent results connected to the uniqueness problem. Recall that the Cauchy problem for the Navier-Stokes equations on $\mathbb{R}^n$ is
\begin{equation} \label{NS_Rn}
\begin{cases}
\partial_t u + u\cdot\nabla u - \nu\Delta u + \nabla p = 0, & \\
\mathrm{div}~u = 0, & \\
u(0, x) = u_0(x),
\end{cases}
\end{equation}
where $u : \mathbb{R} \times \mathbb{R}^n \rightarrow \mathbb{R}^n$ is the velocity vector field, $\nu > 0$ is a viscosity parameter and $p : \mathbb{R} \times \mathbb{R}^n \rightarrow \mathbb{R}$ is the pressure. From now on, we consider only the cases $n = 2, 3$, which are of special interest from the standpoint of classical mechanics and have been studied extensively over the past century.

Among the first fundamental contributions to this problem, going back to $1934$, were the works of Leray \cite{Leray} and Hopf \cite{Hopf}, establishing that $(\ref{NS_Rn})$ admits global weak solutions
\begin{equation}\label{LH_1}
u \in L^{\infty}\big([0, \infty), L^2(\mathbb{R}^n)\big) \cap L^2\big( [0, \infty), \dot{H}^1(\mathbb{R}^n)\big)
\end{equation}
satisfying the energy inequality
\begin{equation}\label{LH_2}
\|u(t, \cdot)\|_{L^2}^2 + 4\int\limits_0^t \|\mathrm{Def}(u(s, \cdot))\|_{L^2}^2 ds \leq \|u_0\|_{L^2}^2,
\end{equation}
where $\mathrm{Def}(u) = \frac{1}{2}\big(\nabla u + (\nabla u)^t \big)$ is the so-called \textit{deformation tensor}. Solutions of $(\ref{NS_Rn})$ satisfying $(\ref{LH_1})$ and $(\ref{LH_2})$ are said to be of \textit{Leray-Hopf} type. 

When $n = 2$, it is known that Leray-Hopf solutions are, in fact, smooth and unique, but this is not known to hold for $n = 3$ (see \cite{Lady1}, \cite{Constantin} for these and related results).

In $1976$, Heywood \cite{Heywood} considered $(\ref{NS_Rn})$ on a domain $D \subseteq \mathbb{R}^n$, under suitable boundary conditions, and showed that the uniqueness problem is related to the geometry of $D$. For instance, on
\begin{equation}
D = \{ (x, y, z) \in \mathbb{R}^3 : x \neq 0,~\text{or}~y^2 + z^2 < 1 \}
\end{equation}
there are nonunique solutions, whereas on bounded domains uniqueness always holds.

In view of these developments, it is natural to consider the more general question of uniqueness of solutions of $(\ref{NS_Rn})$ on a general Riemannian manifold $(M, g)$.

In what follows, we are concerned only with the case when $(M, g)$ is not flat. In this context, the system $(\ref{NS_Rn})$ takes the form
\begin{equation} \label{NS_M}
\begin{cases}
\partial_t u + \nabla_uu - \nu\mathcal{L}u + \nabla p = 0, & \\
\mathrm{div}~u = 0, & \\
u(0, x) = u_0(x),
\end{cases}
\end{equation}
where $\nabla$ is the Levi-Civita connection, and $\mathcal{L}$ is an operator generalizing the Euclidean Laplacian (see below). When $\nu = 0$ in the first equation of $(\ref{NS_M})$, we obtain the well-known Euler equations
\begin{equation}\label{euler_eq}
\begin{cases}
\partial_t u + \nabla_uu + \nabla p = 0, & \\
\mathrm{div}~u = 0, & \\
u(0, x) = u_0(x),
\end{cases}
\end{equation}
Throughout the literature, we encounter several natural choices for the operator $\mathcal{L}$. The first one is the Bochner (or \textit{rough}) Laplacian, given by $\Delta_B = \nabla^*\nabla$, where $\nabla^*$ is the formal adjoint of $\nabla$. Another possibility is to take $\mathcal{L}$ to be the Hodge (or \textit{damped}) Laplacian (cf. \cite{ANA}), given by $\Delta_H = dd^* + d^*d$, where $d$ and $d^*$ are the exterior derivative and its formal adjoint. Finally, one can also take (see \cite{EM}) $\Delta_D = 2\mathrm{Def}^*\mathrm{Def}$, where $\mathrm{Def}$ is the deformation tensor.

The operators above, when acting on a divergence-free $1$-form $w$, are related by
\begin{equation}
\Delta_D(w) = \Delta_B(w) - \mathrm{Ric}(w, \cdot) = \Delta_H(w) - 2\mathrm{Ric}(w, \cdot),
\end{equation}
where $\mathrm{Ric}$ is the Ricci tensor and we identify $1$-forms with vector fields via $g$.

According to Ebin and Marsden \cite{EM}, when $(M, g)$ is Einstein (i.e., $\mathrm{Ric} = \lambda g$, for some constant $\lambda$), the appropriate operator is $\mathcal{L} = \Delta_D$. Nevertheless, the other choices also appear in the literature. In particular, the system $(\ref{NS_M})$ with $\mathcal{L} = \Delta_H$ (referred to in \cite{CC} as the  \textit{modified} Navier-Stokes equations), which we consider here, was studied in a number of works (e.g., \cite{AB}, \cite{RW}, \cite{CMT}, \cite{ANA}, \cite{IL1}, \cite{CC}).

The case of compact Riemannian manifolds $(M, g)$, with or without boundary, has been studied by several authors (e.g., \cite{MM}, \cite{IL1}, \cite{IL2}, \cite{EM}). For instance, in \cite{EM} it was proved that $(\ref{NS_M})$ is well-posed on Sobolev spaces, in any dimension.

For non-compact Riemannian manifolds, however, the situation is quite different. In sharp contrast with the classical results in $\mathbb{R}^2$ mentioned above, Chan and Czubak \cite{CC} proved in $2010$ that there are nonunique Leray-Hopf solutions of $(\ref{NS_M})$ on the two-dimensional hyperbolic space $\mathbb{H}^2(-a^2)$. In \cite{CC}, the authors use $\mathcal{L} = \Delta_D$ (since $\mathbb{H}^2(-a^2)$ is Einstein), and also prove a nonuniqueness result where $\mathcal{L} = \Delta_H$ for more general negatively curved two-dimensional manifolds.

The question was raised in \cite{CC} as to whether their methods could be extended to three-dimensional manifolds.

In the case of $\mathbb{H}^3(-a^2)$, the above question was answered negatively by Khesin and Misiolek \cite{KM} in $2012$, who showed that the nonuniqueness phenomenon observed in \cite{CC} is essentially due to the Hodge-Kodaira decomposition, and that similar techniques cannot produce nonunique solutions on $\mathbb{H}^n(-a^2)$, for $n > 2$. Some of their findings can be summarized in the following theorem, which incorporates an earlier result of Dodziuk \cite{Dod}.

\begin{theorem}\label{km_thm}
Given a Riemannian manifold $(M, g)$, every $L^2$ harmonic vector field that belongs to the Sobolev space $W^{1,4}(TM)$ defines a time-independent solution of the Euler equations $(\ref{euler_eq})$. Furthermore, there is an infinite-dimensional space of such vector fields on the Hyperbolic plane $\mathbb{H}^2(-a^2)$. On the other hand, for $n \geq 3$ any $L^2$ harmonic vector field on $\mathbb{H}^n(-a^2)$ vanishes identically.
\end{theorem}

It is worth pointing out here that the space $W^{1,4}(TM)$ does not appear explicitly in \cite{KM} due to a typo in the statement of Proposition $1.2$ of that reference. As one sees directly from its proof, Proposition $1.2$ of \cite{KM} requires precisely this decay.

The relationship between Theorem $\ref{km_thm}$ and the Cauchy problem for the Navier-Stokes equations $(\ref{NS_M})$ is as follows. Let $v = \nabla \phi$ be a time-independent solution of the Euler equations $(\ref{euler_eq})$, for some function $\phi : M \rightarrow \mathbb{R}$. If $v$ decays as in Theorem $\ref{km_thm}$ and has finite dissipation (i.e., finite $\dot{H}^1$ norm), one obtains infinitely many Leray-Hopf solutions $u = f(t)v$ of $(\ref{NS_M})$, for carefully chosen real-valued functions $f(t)$. Because of the finite dissipation condition, one ends up having to require both $W^{1, 2}$ and $W^{1, 4}$ norms to be finite.

This was the basic technique employed in \cite{CC} to obtain a nonuniqueness result in $\mathbb{H}^2(-a^2)$. In particular, one sees from Theorem $\ref{km_thm}$ that the methods of \cite{CC} do not yield nonuniqueness of $(\ref{NS_M})$ in the case of $\mathbb{H}^3(-a^2)$.

Nevertheless, we show that the approach introduced in \cite{CC} can be implemented on certain three-dimensional manifolds of non-constant negative curvature.

As a matter of fact, the manifolds introduced by Anderson \cite{And} in his counterexample to the Dodziuk-Singer conjecture can be shown, after some delicate estimates, to carry nontrivial harmonic vector fields with sufficient decay:

\begin{theorem}\label{thm_1}
The Anderson manifold $(M^3, g)$ has the property that the space of harmonic vector fields in $W^{1, 2}(TM) \cap W^{1,4}(TM)$ is infinite-dimensional.
\end{theorem}

We emphasize here that, in order to prove Theorem $\ref{thm_1}$, it is not enough to repeat the construction of the harmonic forms as in \cite{And}. We first need  to produce a new, nontrivial supersolution to a certain PDE (system $(\ref{bvp})$). This is done in Lemma $\ref{lem:supersolution}$, and provides stronger decay properties for the solutions of the given PDE. After that, we show that these decay properties are sharp (Lemma $\ref{lem_psi}$), in the sense that they are exactly what is needed for the harmonic forms to lie in the space $W^{1, 2}(T^*M) \cap W^{1, 4}(T^*M)$. Finally, we carry out the Sobolev estimates in full detail.

As an immediate corollary of Theorems $\ref{km_thm}$ and $\ref{thm_1}$, and the construction in \cite{CC} mentioned above, we obtain the following nonuniqueness result for the Navier-Stokes equations.

\begin{theorem}\label{thm_2}
The Cauchy problem $(\ref{NS_M})$, with $\mathcal{L} = \Delta_H$, is ill-posed on $(M^3, g)$ in the sense that there are smooth initial conditions $u_0$ for which one has infinitely many smooth Leray-Hopf solutions of $(\ref{NS_M})$.
\end{theorem}

One could state a version of Theorem $\ref{thm_2}$ for the Euler equations $(\ref{euler_eq})$. However, our solutions do not conserve energy (i.e., the $L^2$ norm is not constant), and this is usually regarded as an essential feature of solutions of the Euler equations.

The paper is organized as follows. In Section $2$, besides reviewing the construction of the manifold $(M^3, g)$, we prove the existence of a new supersolution that is more refined than the one used in \cite{And} to create nontrivial $L^2$ harmonic $1$-forms.

The proof of Theorem $\ref{thm_1}$ is completed in Section $3$, where we carry out the estimates necessary to prove that the $1$-forms constructed in Section $2$ satisfy the hypothesis of Theorem $\ref{km_thm}$.

For other results concerning the Navier-Stokes equations on a Riemannian manifold, see \cite{CC} and references therein.

\subsection*{Acknowledgements}
The author wishes to thank Prof. Gerard Misiolek for suggesting the problem and for many inspiring discussions, and Prof. Frederico Xavier for his interest and helpful comments on this work. This research was supported by the Richard and Peggy Notebaert Fellowship.

\section{The construction of $(M^3, g)$}
\label{sec:2} 

Throughout this paper, whenever $E \rightarrow M$ is a vector bundle equipped with a Riemannian metric, we denote by $L^p(E)$ the space of $p$-integrable sections of $E$ and by $W^{s, p}(E)$ the $L^p$-based Sobolev space of sections of $E$, with smoothness index $s$. See \cite{Liviu} for more details.

In this section, we recall Anderson's construction of a complete Riemannian three-dimensional manifold $(M^3, g)$ and perform the modifications needed in our case. The latter are necessary in order to ensure that the constructed harmonic vector fields obey appropriate decay conditions at infinity. Our basic reference for this section is \cite{And}.


Consider the set
\begin{gather}\begin{split}
H = \{ (x, y) \in \mathbb{C} : x^2 + y^2 < 1,~x > 0 \}.
\end{split}\end{gather}
Fix a number $a > 1$ and equip $H$ with the hyperbolic metric $ds_a^2$ of sectional curvature $-a^2$. In cartesian coordinates $z = x + iy$, it takes the form
\begin{equation}
ds_a^2 = 4\frac{dx^2 + dy^2}{a^2(1 - |z|^2)^2}.
\end{equation}
Consider a warped product $H \times_f S^1$, where $S^1$ is the circle and $f : H \rightarrow \mathbb{R}$ is a smooth function to be defined. Let $L$ be the line segment from $(0, -1)$ to $(0, 1)$ and let $\rho_a(z)$ be the distance, with respect to $ds_a^2$, from the point $z \in H$ to $L$. More explicitly, we have (see \cite[p.~162]{Be})
\begin{equation}
\rho_a(z) = \frac{1}{a} \mathrm{arcsinh} \left( \frac{2x}{1 - |z|^2} \right).
\end{equation}
Define $f$ by
\begin{equation} \label{formula:f}
f(z) = \mathrm{sinh} \circ \left(\rho_a(z)\right).
\end{equation}
The explicit formula for $f$ will be used later on, but for now we only remark that from the above expression it follows that $f$ extends smoothly to $H \cup L$, and that the corresponding warped metric extends to $M = H \times_f S^1$, which is topologically just $\mathbb{R}^3$ (for the metric, see \cite[p.~12]{Pete}). Using cartesian coordinates $z = x + iy$ on $H$ and the standard angular coordinate $\theta$ on $S^1$, the metric on $M$ reads
\begin{equation}\label{metric_M}
ds_M^2 = 4\frac{dx^2 + dy^2}{a^2(1 - |z|^2)^2} + f^2d\theta^2.
\end{equation}


To obtain a large space of harmonic $1$-forms on $M$, we first decompose the Hodge Laplacian $\Delta_M$ acting on $1$-forms on $M$ as
\begin{equation}\label{eq:laplacian_M}
\Delta_M = \Delta_{\mathbb{H}^2\text{\tiny $(-a^2)$}} - d \circ i_F + i_F \circ d,
\end{equation}
where $F = \nabla \log (f)$, $i_F$ denotes interior multiplication by $F$ and $\nabla$ is the gradient of the Riemannian metric on $M$. Since $f = 0$ on $L$, the above formula only makes sense on $H \times_f S^1$. However, if we require that solutions of the equation
\begin{equation}\label{eq:laplacian_M_0}
\Delta_{\mathbb{H}^2\text{\tiny $(-a^2)$}}(\omega) - d(\omega(F)) + (d\omega)(F, \cdot) = 0,
\end{equation}
satisfy a suitable Neumann condition at $L$, then they can be extended smoothly to $M$, and the extension will be harmonic. First, we require $\omega$ to be invariant under isometric reflection through $L$. Second, we require that $\omega$ be invariant under the isometric action of the $S^1$-factor. This reduces the problem to solving equation $(\ref{eq:laplacian_M_0})$ on $H$.
\newline\newline \textit{Remark}. From this point on, we specialize to the case $a = 2$ for simplicity. Other values of $a > 1$ would work as well, with minor modifications.
\newline\newline Assume now that $\omega = du$. Let $\Delta = \Delta_{\mathbb{H}^2\text{\tiny $(-4)$}}$. Then $(\ref{eq:laplacian_M_0})$ becomes
\begin{equation}\label{eq:laplacian_M_1}
d\left(\Delta u - du(F)\right) = 0,
\end{equation}
so that it is enough to solve
\begin{equation}\label{eq:laplacian_M_2}
\Delta u - du(F) = 0
\end{equation}
on $H$. We would like to rewrite this equation on the strip
\begin{equation}
\Omega = \{ (r, s) \in \mathbb{C} : 0 < s < \pi/2 \},
\end{equation}
where it takes a simpler form. Thus, we introduce the biholomorphisms $\varphi : \Omega \rightarrow H$ and $\varphi^{-1} = \psi : H \rightarrow \Omega$ given by
\begin{gather}\begin{split}\label{eq_phipsi}
\varphi(w) = i\dfrac{e^{-w} - 1}{e^{-w} + 1},~~\psi(z) = \log\left(\dfrac{i - z}{i + z}\right).
\end{split}\end{gather}
Note that $\psi(L) = \mathbb{R} \times \{ 0 \}$, and that the image of the full unit disk under $\psi$ is
\begin{equation}
\widetilde{\Omega} = \{ (r, s) : -\pi/2 < s < \pi/2 \}.
\end{equation}
Letting $v = u \circ \varphi$, $(\ref{eq:laplacian_M_2})$ can be rewritten on $\Omega$ as
\begin{equation}\label{eq:laplacian_M_3}
L_{\Omega}(v) = \frac{\partial^2 v}{\partial r^2} + \frac{\partial^2 v}{\partial s^2} + \dfrac{\partial \big(\log(\widetilde{f})\big)}{\partial s}\frac{\partial v}{\partial s} = 0,
\end{equation}
where $\widetilde{f} = f \circ \varphi$. This means that $u$ is a solution of $(\ref{eq:laplacian_M_2})$ if and only if $v$ is a solution of $(\ref{eq:laplacian_M_3})$. A long but straightforward computation using $(\ref{formula:f})$ gives
\begin{gather}\begin{split}\label{formula:f_on_strip}
\widetilde{f}(r, s) = \frac{1}{2}\dfrac{(1 + \sin(s))^{1/2} - (1 - \sin(s))^{1/2}}{\cos(s)^{1/2}}.
\end{split}\end{gather}
Observe that $\widetilde{f}$ depends only on one variable. Each solution of $(\ref{eq:laplacian_M_3})$ gives a harmonic $1$-form on $H$. 
Consider the mixed boundary value problem
\begin{gather}\begin{split} \label{bvp}
\begin{cases}
L_{\Omega}(v) = 0, & \\
\dfrac{\partial v}{\partial s}(r, 0) = 0, & \\
v(r, \pi/2) = v_0(r),
\end{cases}
\end{split}\end{gather}
where $v_0 : \mathbb{R} \rightarrow \mathbb{R}$ is a smooth function with uniformly bounded derivatives of all orders. It was shown in \cite{And} that $(\ref{bvp})$ has a unique solution $v$ which can be extended by reflection to
\begin{equation}
\widetilde{\Omega} = \{ (r, s) : -\pi/2 < s < \pi/2 \}.
\end{equation}
This solution is smooth in $\widetilde{\Omega}$ and continuous up to the boundary
\begin{equation}
\partial{}\widetilde{\Omega} = \{ (r, s) : s = \pi/2 \} \cup \{ (r, s) : s = -\pi/2 \}.
\end{equation}
In particular, by construction, $v = v_0$ on the boundary of $\widetilde{\Omega}$. By the maximum principle applied to the operator $L_{\Omega}$ (see \cite[p.~35]{GT}), constructing appropriate supersolutions of $(\ref{bvp})$ yields bounds on $v$ and its first and second derivatives. More precisely, if $v^+$ is a positive supersolution of $(\ref{bvp})$ and $v$ is chosen such that $v(r, \pi/2) \leq v^+(r, \pi/2)$, then
\begin{equation}\label{bounds_v_pre}
|v| + |v_r| + |v_s| + |v_{rr}| + |v_{rs}| + |v_{ss}| \leq \text{\it const.}~v^+
\end{equation}
holds \textit{pointwise}. It should be remarked that $(\ref{bounds_v_pre})$ relies on estimates obtained from the method used to prove the existence of solutions to $(\ref{bvp})$ in \cite{And}, rather than a general argument.
\newline Let
\begin{equation}
v^+(r, s) = e^{-\delta |r|}p(s),
\end{equation}
where $p : [0, \pi/2] \rightarrow (0, +\infty)$ is a smooth function and $\delta > 0$. We can see directly from $(\ref{eq:laplacian_M_3})$ that $v^+$ is a positive supersolution of $(\ref{bvp})$ (with initial condition $v^+_0(r) = e^{-\delta |r|}p(\pi/2)$) if and only if
\begin{equation} \label{supersolution}
\begin{cases}
p''(s) + \left(\dfrac{~\widetilde{f}_s}{\widetilde{f}}\right)p'(s) + \delta^2 p(s) \leq 0, & \\
p'(0) = 0.
\end{cases}
\end{equation}
In \cite{And}, the author uses $p(s) = -c_1s^2 + c_2$, for carefully chosen $c_1, c_2$ and $\delta > 0$ sufficiently small. In that case, for optimal $c_1$ and $c_2$, the number $\delta$ has to be less than $\sqrt{6}/\pi \approx 0.78$. However, our construction requires $\delta$ to be at least $1$, as shown by Lemma $\ref{lem_psi}$.  For this, we need the following
\begin{lem}\label{lem:supersolution}
The function
\begin{equation}
v^+(r, s) = e^{-|r|}(8s^4 - 50s^2 + 75),
\end{equation}
is a supersolution of $(\ref{bvp})$.
\end{lem}
\noindent \textit{Proof}. We must check that $p(s) = 8s^4 - 50s^2 + 75$ is positive on $[0, \pi/2]$ and that it satisfies $(\ref{supersolution})$ with $\delta = 1$. Positivity follows from the fact that the roots of $p(s)$ are $\pm \sqrt{5/2}, \pm \sqrt{15}/2$, and $p(0) > 0$. The inequality in $(\ref{supersolution})$ becomes
\begin{equation}\label{ineq_lemma}
(8s^4 + 46s^2 - 25) + G(s)(32s^3 - 100s) \leq 0,
\end{equation}
where
\begin{equation}\label{eq_G}
G(s) = \widetilde{f}_s / \widetilde{f} = \frac{1}{2}\left(\frac{\sqrt{1 + \sin(s)} + \sqrt{1 - \sin(s)}}{\sqrt{1 + \sin(s)} - \sqrt{1 - \sin(s)}} + \tan(s)\right).
\end{equation}
It is clear that $G(s) > 0$ and $p_1(s) = 32s^3 - 100s < 0$ on $(0, \pi/2)$. Only $p_2(s) = 8s^4 + 46s^2 - 25$ changes sign, so we break the inequality $(\ref{ineq_lemma})$ in two cases:
\begin{equation*}
\begin{cases}
0 < s \leq \sqrt{2}/2 & \\
\sqrt{2}/2 < s < \pi/2
\end{cases}
\end{equation*}
In the first case, $p_2(s)$ is negative, so $(\ref{ineq_lemma})$ obviously holds. In the second case, we have $1/2 < \sin(s) \leq 1$, and this gives
\begin{equation}\label{ineq_2}
G(s) \geq \frac{1}{2}\left(1 + \frac{1}{2\cos(s)}\right),
\end{equation}
so that
\begin{gather}\begin{split}\label{ineq_3}
(8s^4 + 46s^2 - 25) + G(s)(32s^3 - 100s) &\leq (8s^4 + 46s^2 - 25) \\
&+(16s^3 - 50s) + \frac{8s^3 - 25s}{\cos(s)},
\end{split}\end{gather}
and this last expression is negative on the interval $(\sqrt{2}/2, \pi/2)$. In either case, $(\ref{ineq_lemma})$ holds, thus proving the Lemma. $\square$
\newline\newline We will only work with solutions $v$ of $(\ref{bvp})$ obeying $(\ref{bounds_v_pre})$. Lemma $\ref{lem:supersolution}$ implies the pointwise estimate
\begin{equation}\label{bounds_v}
|v| + |v_r| + |v_s| + |v_{rr}| + |v_{rs}| + |v_{ss}| \leq \text{\it const.}~e^{-|r|}.
\end{equation}
We finish this section with two Lemmas. The last one is the main ingredient needed in the next section to establish that $du \in W^{1, 2}(T^*M) \cap W^{1, 4}(T^*M)$.

\begin{lem}\label{lem_psi}
Let $z \in H$ and $\psi(z) = \psi_1(z) + i\psi_2(z)$, where $\psi_1, \psi_2$ are real-valued and $\psi$ is the function defined in $(\ref{eq_phipsi})$. The functions
\begin{equation}
f_1(z) = (1 - |z|^2)|\psi'(z)|\hspace{2mm}\text{and}\hspace{2mm}f_2(z) = e^{-\delta|\psi_1(z)|}|\psi'(z)|
\end{equation}
are bounded on $H$, provided that $\delta$ is at least $1$. If $\delta < 1$, then $f_2$ is not bounded.
\end{lem}
\noindent \textit{Proof}. That $f_1$ is bounded follows directly from the triangle inequality. As for $f_2$, note that $|\psi'(z)|$ has only two singularities,  located at $\pm i$. We analyze the one at $z = i$, the other one being completely analogous. We have
\begin{equation}
e^{-2\delta|\psi_1(z)|}|\psi'(z)|^2 = \text{\it const.}~\dfrac{(x^2 + (1-y)^2)^{\delta}}{(1 + x^2 - y^2)^2 + 4x^2y^2},
\end{equation}
and switching to polar coordinates $(r, \theta)$,
\begin{gather}\begin{split}\label{eq_delta_1}
e^{-2\delta|\psi_1(z)|}|\psi'(z)|^2 &= \text{\it const.}~\dfrac{(1 - 2r\sin(\theta) + r^2)^{\delta}}{1 + 2r^2\cos(2\theta) + r^4} \approx \text{\it const.}~\dfrac{(1-r)^{2\delta}}{(1+r)^2(1-r)^2} \\
&\approx \text{\it const.}~\dfrac{1}{(1-r)^{2-2\delta}}.
\end{split}\end{gather}
since $\theta \approx \pi/2$, $r \approx 1$. ~$\square$

\begin{lem}\label{lem_bounds_u}
For $u = v \circ \psi$, where $v$ is a solution of $(\ref{bvp})$ satisfying $(\ref{bounds_v})$, we have the pointwise estimates
\begin{equation}\label{bounds_u_1}
|u_x| + |u_y| \leq \text{\it const.}
\end{equation}
and
\begin{equation}\label{bounds_u_2}
|u_{xx}| + |u_{xy}| + |u_{yy}| \leq \text{\it const.}~|\psi'|.
\end{equation}
\end{lem}

\noindent \textit{Proof}. From Lemma $\ref{lem_psi}$ and $(\ref{bounds_v})$, we have
\begin{equation}
|u_x| + |u_y| \leq \text{\it const.}~e^{-|\psi_1|}|\psi'| \leq \text{\it const.}
\end{equation}
For the second derivatives of $u$, using the Cauchy-Riemann equations and the bound $|\psi''(z)| \leq |\psi'(z)|^2$, 
\begin{gather}\begin{split}
|u_{xx}| &= |v_{rr}(\psi_1)_x^2 + 2v_{rs}(\psi_1)_x(\psi_2)_x + v_{ss}(\psi_2)_x^2 + v_r(\psi_1)_{xx} + v_{s}(\psi_2)_{xx}| \\
&\leq \mathrm{\it \textit const.}~e^{-{|\psi^1}|}~\big((\psi^1_x)^2 + 2|\psi^1_x\psi^2_x| + (\psi^2_x)^2 + |\psi^1_{xx}| + |\psi^2_{xx}|         \big) \\
&\leq \mathrm{\it \textit const.}~e^{-{|\psi^1}|}~\big(|\psi'|^2 + |\psi''|\big) \\
&\leq \mathrm{\it \textit const.}~e^{-{|\psi^1}|}~|\psi'|^2 \\
&\leq \mathrm{\it \textit const.}~|\psi'|.
\end{split}\end{gather}
The same estimate holds for $|u_{xy}|$ and $|u_{yy}|$.~$\square$

\section{Sobolev estimates} 
\label{sec:3} 

In this section, we show that the $W^{1,2}$ and $W^{1, 4}$ norms of $du$ are finite. Recall that, up to a set of measure zero, $M = H \times_f S^1$. Since $S^1$ is compact and the integrands will not depend on $\theta$, we can reduce all computations to $H$.
\newline\newline Denote by $\|~\|_g$ the norm with respect to the metric $(g_{ij})$ of $M$, so that from $(\ref{metric_M})$ we have
\begin{equation*}
\|dx\|_g = \|dy\|_g = (1 - |z|^2)\hspace{2.5mm}\text{and}\hspace{2.5mm}\|d\theta\|_g = 1/f,
\end{equation*}
where $z = x +iy$ as before. Then, using $(\ref{bounds_u_1})$,

\begin{gather}\begin{split} \label{l2_est}
\int\limits_H \|du\|_g^2 f \frac{dx~dy}{(1 - |z|^2)^2} &= \int\limits_H (u_x^2 + u_y^2) f ~dx~dy \leq \text{\it const.}~\int\limits_H f ~dx~dy.
\end{split}\end{gather}

\subsubsection*{$L^2$ estimate}
We show explicitly that the last integral in $(\ref{l2_est})$ is finite. Expanding $(\ref{formula:f})$ gives
\begin{gather*}\begin{split}
f(x, y) = -\frac{1}{2}\left( \frac{2x}{1 - |z|^2} + \sqrt{1 + \dfrac{4x^2}{(1 - |z|^2)^2}} \right)^{-1/2} +\frac{1}{2}\left( \frac{2x}{1 - |z|^2} + \sqrt{1 + \dfrac{4x^2}{(1 - |z|^2)^2}} \right)^{1/2}.
\end{split}\end{gather*}
The first term above is bounded by $2$ in absolute value, since
\begin{gather}\begin{split}
\frac{2x}{1 - |z|^2} + \sqrt{1 + \dfrac{4x^2}{(1 - |z|^2)^2}} \geq 1.
\end{split}\end{gather}
As for the second term, note that
\begin{gather}\begin{split}
\left( \frac{2x}{1 - |z|^2} + \sqrt{1 + \dfrac{4x^2}{(1 - |z|^2)^2}} \right)^{1/2} &\leq \left( \frac{4x}{1 - |z|^2} + 1 \right)^{1/2} \\
&\leq \sqrt{\dfrac{4x}{1 - |z|^2}} + 1.
\end{split}\end{gather}
Using polar coordinates,
\begin{gather}\begin{split}\label{eq:l2}
\int\limits_H f~dx~dy \leq \text{\it const.} + \text{\it const.}~\int_{-\pi/2}^{\pi/2} \int_0^1 \frac{2r}{\sqrt{1 - r^2}} dr d\theta < \infty.
\end{split}\end{gather}
Therefore, we have $du \in L^2(T^*M)$.

\subsubsection*{$W^{1, 2}$ and $W^{1, 4}$ estimates}
From $(\ref{eq:l2})$ and the fact that
\begin{gather}\begin{split}
\|du\|_g^4 = (1 - |z|^2)^4(u_x^2 + u_y^2)^2 \leq \text{\it const.}~(1 - |z|^2)^2,
\end{split}\end{gather}
we see immediately that $du \in L^4(T^*M)$. It remains to show that
\begin{equation*}
\nabla du \in L^2(T^*M \otimes T^*M) \cap L^4(T^*M \otimes T^*M),
\end{equation*}
where $\nabla$ denotes the Levi-Civita connection acting on $1$-forms. First, we compute and find an upper bound for $\|\nabla du\|_g^2 = \langle \nabla du, \nabla du \rangle$. To do this, we break up the expression into three terms:
\begin{gather*}\begin{split}
\|\nabla du\|_g^2 &= \| u_{xx}dx \otimes dx + u_{xy} dx \otimes dy + u_{xy} dy \otimes dx + u_{yy} dy \otimes dy \|_g^2 \\
&+2\Big\langle u_{xx}dx \otimes dx + u_{xy} dx \otimes dy + u_{xy} dy \otimes dx + u_{yy} dy \otimes dy,~ u_x\nabla dx + u_y \nabla dy \Big\rangle \\
&+\|u_x\nabla dx + u_y \nabla dy \|_g^2
\end{split}\end{gather*}
and then estimate each one separately.
\newline\newline Let us start from the last term containing only the lower order derivatives of $u$ appearing in $\nabla du$. Recall that, using Einstein's summation convention, we have
\begin{equation*}
\nabla dx^k = -\Gamma^k_{ij} dx^i \otimes dx^j,
\end{equation*}
where $\Gamma^k_{ij}$ are the Christoffel symbols, given in terms of the metric $(g_{ij})$ of $M$ by the well-known formula
\begin{equation*}
\Gamma^k_{ij} = \frac{1}{2} g^{kl} \left( \frac{\partial g_{il}}{\partial x_j} + \frac{\partial g_{jl}}{\partial x_i} - \frac{\partial g_{ij}}{\partial x_l} \right).
\end{equation*}
From this we compute
\begin{gather}\begin{split}
\nabla dx &= \frac{2x}{1 - |z|^2} dx \otimes dx + \frac{2y}{1 - |z|^2} dx \otimes dy + \frac{2y}{1 - |z|^2} dy \otimes dx \\
&+ \frac{-2x}{1 - |z|^2} dy \otimes dy -(1 - |z|^2)^2ff_x d\theta \otimes d\theta, \\ \\
\nabla dy &= \frac{-2y}{1 - |z|^2} dx \otimes dx + \frac{2x}{1 - |z|^2} dx \otimes dy +\frac{2x}{1 - |z|^2} dy \otimes dx \\
&+ \frac{2y}{1 - |z|^2} dy \otimes dy - (1 - |z|^2)^2ff_y d\theta \otimes d\theta.
\end{split}\end{gather}
To obtain an upper bound for $\| u_x \nabla dx + u_y \nabla dy \|_g^2$, we split
\begin{gather*}\begin{split}
\| u_x \nabla dx + u_y \nabla dy \|_g^2 = \| u_x \nabla dx \|_g^2 + \| u_y \nabla dy \|_g^2 + 2 \langle u_x \nabla dx, u_y \nabla dy \rangle
\end{split}\end{gather*}
further into three parts:
\begin{gather*}\begin{split}
\| u_x \nabla dx \|_g^2 &= u_x^2(1 - |z|^2)^2\left(4x^2 +4y^2 + 4y^2 + 4x^2 + (1 - |z|^2)^2\frac{f_x^2}{f^2}\right), \\
\| u_y \nabla dy \|_g^2 &= u_y^2(1 - |z|^2)^2\left(4y^2 +4x^2 + 4x^2 + 4y^2 + (1 - |z|^2)^2\frac{f_y^2}{f^2}\right), \\
2\langle u_x \nabla dx, u_y \nabla dy \rangle &= 2u_xu_y(1 - |z|^2)^2\left( -4xy + 4xy + 4xy - 4xy + (1 - |z|^2)^2\frac{f_xf_y}{f^2} \right) \\
&= 2u_xu_y(1 - |z|^2)^4\frac{f_xf_y}{f^2}.
\end{split}\end{gather*}
Adding all three terms gives us
\begin{gather}\begin{split} \label{eq_nabla_lowerorder}
\|u_x \nabla dx + u_y \nabla dy \|_g^2 &= 8|z|^2(u_x^2 + u_y^2)(1 - |z|^2)^2 \\
&+\frac{(1 - |z|^2)^4}{f^2}\big(u_xf_x + u_yf_y\big)^2.
\end{split}\end{gather}
Now, recall that $u$ satisfies $(\ref{eq:laplacian_M_2})$, that is, $\Delta u = du(F)$, where $F = \nabla \log(f)$. In coordinates, this becomes
\begin{gather}\begin{split}
\Delta u = du(\nabla \log(f)) = \frac{g^{11}}{f}u_xf_x + \frac{g^{22}}{f}u_yf_y = \frac{(1 - |z|^2)^2}{f}(u_xf_x + u_yf_y),
\end{split}\end{gather}
where $(g^{ij})$ are the components of the inverse of $g$, so that
\begin{gather}\begin{split}
u_xf_x + u_yf_y = \frac{f}{(1 - |z|^2)^2}\Delta u.
\end{split}\end{gather}
\normalsize
Substituting this back into $(\ref{eq_nabla_lowerorder})$ gives
\begin{gather}\begin{split} \label{eq_nabla_lowerorder2}
\|u_x \nabla dx + u_y \nabla dy \|_g^2 &= 8|z|^2(u_x^2 + u_y^2)(1 - |z|^2)^2 + (\Delta u)^2.
\end{split}\end{gather}
Thus, it remains to estimate $|\Delta u|$. Using $(\ref{bounds_u_2})$ and Lemma $\ref{lem_psi}$, we get
\begin{gather}\begin{split} \label{eq_laplacian_u}
|\Delta u| &= |g^{11}u_{xx} + g^{22}u_{yy}| \\
&\leq \mathrm{\it \textit const.}~(1 - |z|^2)^2|\psi'(z)| \\
&\leq \mathrm{\it \textit const.}~(1 - |z|^2).
\end{split}\end{gather}
Going back to $(\ref{eq_nabla_lowerorder2})$ and using $(\ref{bounds_u_1})$, we get
\begin{gather}\begin{split} \label{eq_nabla_lowerorder3}
\|u_x \nabla dx + u_y \nabla dy \|_g^2 \leq \text{\it const.}~(1 - |z|^2)^2. \\
\end{split}\end{gather}
This takes care of the term containing only first order derivatives of $u$ in $\langle \nabla du, \nabla du \rangle$. For the term containing only second order derivatives, using Lemma $\ref{lem_psi}$ and $(\ref{bounds_u_2})$ gives
\begin{gather}\begin{split}
&\| u_{xx}dx \otimes dx + u_{xy} dx \otimes dy + u_{xy} dy \otimes dx + u_{yy} dy \otimes dy \|_g^2 = \\
&=(u_{xx}^2 + 2u_{xy}^2 + u_{yy}^2)(1 - |z|^2)^4 \\
&\leq \mathrm{\it const.}~|\psi'(z)|^2(1 - |z|^2)^4 \\
&\leq \mathrm{\it const.}~(1 - |z|^2)^2,
\end{split}\end{gather}
which is the same bound as before. Finally, for the mixed term that contains first and second derivatives of $u$ in $\langle \nabla du, \nabla du \rangle$, we again use $(\ref{bounds_u_1})$, $(\ref{bounds_u_2})$ and Lemma $\ref{lem_psi}$:
\begin{gather*}\begin{split}
&2\Big\langle u_{xx}dx \otimes dx + u_{xy} dx \otimes dy + u_{xy} dy \otimes dx + u_{yy} dy \otimes dy, u_x\nabla dx + u_y \nabla dy \Big\rangle \\
&\leq 2(1-|z|^2)^3\Big[|u_x|\Big(2|u_{xx}| + 4|u_{xy}| + 2|u_{yy}|\Big)  +  |u_y|\Big(2|u_{xx}| + 4|u_{xy}| + 2|u_{yy}|\Big)\Big] \\
&\leq \mathrm{\it const.}~(1-|z|^2)^3\Big[|u_x||\psi'(z)|  +  |u_y||\psi'(z)|\Big] \\
&\leq \mathrm{\it const.}~(1-|z|^2)^3|\psi'(z)| \\
&\leq \mathrm{\it const.}~(1-|z|^2)^2.
\end{split}\end{gather*}
Combining all the above estimates and squaring leads to
\begin{gather}\begin{split}
\|\nabla du\|_g^4 \leq \mathrm{\it const.}~(1 - |z|^2)^4,
\end{split}\end{gather}
so that
\begin{gather}\begin{split}
\int\limits_H \|\nabla du\|_g^4~f~\frac{dx~dy}{(1 - |z|^2)^2} \leq \mathrm{\it const.}~\int\limits_H f~dx~dy
\end{split}\end{gather}
and this integral is finite (see $(\ref{eq:l2})$). This shows that $du \in W^{1, 4}(T^*M)$.
\newline Moreover, one sees directly from the above estimate that
\begin{gather}\begin{split}
\int\limits_H \|\nabla du\|_g^2~f~\frac{dx~dy}{(1 - |z|^2)^2} < \infty,
\end{split}\end{gather}
so that finite dissipation, i.e., $du \in \dot{H}^1(T^*M)$, is also satisfied. Combined with the $L^2$ estimate, we have $du \in W^{1,2}(T^*M) \cap W^{1,4}(T^*M)$.

\bibliographystyle{amsplain}

\end{document}